\magnification=1200
\hfuzz=3pt
\overfullrule=0mm



\font\titrefont=cmbx10 at 15pt
\font\sectionfont=cmbx10 scaled\magstep1

\font\notefont=cmr7



\font\tensymb=msam9
\font\fivesymb=msam5 at 5pt
\font\sevensymb=msam7  at 7pt
\newfam\symbfam
\scriptscriptfont\symbfam=\fivesymb
\textfont\symbfam=\tensymb
\scriptfont\symbfam=\sevensymb

\font\sc=cmcsc10 \rm


\def\hfl#1#2{\smash{\mathop{\hbox to 10mm{\rightarrowfill}}
\limits^{\scriptstyle#1}_{\scriptstyle#2}}}

\def\rfl#1#2{\smash{\mathop{\hbox to 10mm{\leftarrowfill}}
\limits^{\scriptstyle#1}_{\scriptstyle#2}}}

\def\vfl#1#2{\llap{$\scriptstyle #1$}\left\downarrow
\vbox to 4mm{}\right.\rlap{$\scriptstyle #2$}}

\def\ufl#1#2{\llap{$\scriptstyle #1$}\left\uparrow
\vbox to 4mm{}\right.\rlap{$\scriptstyle #2$}}

\def\mapright#1{\smash{\mathop{ \longrightarrow}\limits^{#1}}}%

\def\op{{\rm o}}
\def\abs{{\rm a}}
\def\id{{\rm id}}
\def\Ext{{\rm Ext}}
\def\Hom{{\rm Hom}}
\def\Ker{{\rm Ker}}
\def\Im{{\rm Im}}
\def\Tor{{\rm Tor}}
\def\mod{{\rm mod}}

\def\cA{{\cal A}}
\def\cD{{\cal D}}

\def\cI{{\cal I}}
\def\cP{{\cal P}}
\def\cS{{\cal S}}

\def\NN{{\bf N}}
\def\ZZ{{\bf Z}}

\def\ot{{\otimes}}
\def\del{\partial}
\def\eps{\varepsilon}
\def\et{\quad\hbox{et}\quad}

\def\Dem{\noindent {\sc D\'emonstration.--- }}

\def\cqfd{ {\sevensymb {\char 3}}}

\null

\vskip 20pt

\centerline{\titrefont 
Alg\`ebre homologique des N-complexes et}
\medskip
\centerline{\titrefont homologie de Hochschild aux racines de l'unit\'e}

\vskip 30pt

\centerline{\sc Christian Kassel et Marc Wambst}

\vskip 15pt
\noindent
\centerline{\it Institut de Recherche Math\'ematique Avanc\'ee,
Universit\'e Louis Pasteur - C.N.R.S.,}
\centerline{\it 7 rue Ren\'e Des\-cartes, 67084 Strasbourg Cedex, France}

\bigskip
\bigskip
\bigskip
\bigskip

\noindent
{\sc Abstract.}
{\it We set up a homological algebra for $N$-complexes,
which are graded modules together with a degree~$-1$ endomorphism~$d$
satisfying~$d^N=0$. We define $\Tor$- and $\Ext$-groups for $N$-complexes
and we compute them in terms of their classical counterparts ($N=2$).
As an application, we get an alternative definition of the Hochschild
homology of an associative algebra out of an $N$-complex 
whose differential is based on a primitive $N$-th root of unity.}

\bigskip
\noindent
{\sc Mathematics Subject Classification (1991):}
18G25, 18G30, 18G35, 18G99, 05A30, 81R50

\bigskip
\noindent
{\sc Key Words:} 
{\it homological algebra, Hochschild homology, simplicial module,
$q$-calculus}

\bigskip\bigskip\bigskip
\noindent
{\sc R\'esum\'e.}
{\it Nous d\'eveloppons une alg\`ebre homologique pour les $N$-complexes,
c'est-\`a-dire pour des modules gradu\'es munis d'un endomorphisme~$d$
de degr\'e~$-1$ tel que~$d^N=0$. 
Dans ce cadre nous d\'efinissons des groupes $\Tor$ et $\Ext$ 
que nous calculons en fonction des groupes $\Tor$ et $\Ext$ classiques ($N=2$).
Comme application, nous obtenons 
l'homologie de Hochschild d'une alg\`ebre associative
comme l'homologie d'un $N$-complexe dont la diff\'eren\-tielle s'exprime \`a
l'aide d'une racine primitive $N$-i\`eme de l'unit\'e.}

\bigskip
\noindent
{\sc Mots-Cl\'es~:} 
{\it alg\`ebre homologique, homologie de Hochschild, module simplicial,
op\'erateur aux $q$-diff\'erences}

\vfill\eject

\bigskip

\noindent
Le point de d\'epart de ce travail est l'observation que l'homologie de Hochschild
d'une alg\`ebre associative~$A$ au-dessus d'un anneau commutatif~$k$ peut se d\'efinir
\`a partir de diff\'erentielles dans lesquelles $-1$ est remplac\'e par
n'importe quelle racine de l'unit\'e.

Pour tout $n\in \ZZ$, posons $C_n(A) = A^{\ot (n+1)}$ si~$n\geq 0$ et
$C_n(A) = 0$ si~$n<0$.
D\'efinissons $b:C_n(A) \to C_{n-1}(A)$ par
$$b(a_0\ot a_1\ot \cdots \ot a_n)
= \sum_{i=0}^{n-1}\, q^i\, a_0\ot \cdots \ot a_ia_{i+1}\ot \cdots \ot a_n
+ q^n\, a_n a_0\ot a_1\ot \cdots \ot  a_{n-1}
\eqno (0.1)$$
o\`u $q$ est un scalaire et $a_0, a_1, \ldots, a_n \in A$.
Fixons maintenant un entier~$N\geq 2$.
On v\'erifie que~$b^N = 0$ dans les deux cas suivants~:

(i) $q$ est une racine $N$-i\`eme de~$1$ diff\'erente de~$1$~;

(ii)~$q=1$ et $N = 0$ dans l'alg\`ebre~$A$. 

Lorsque $q=-1$ et $N=2$, l'application~$b$
est le bord de Hochschild standard et l'homologie du complexe $(C(A),b)$
est l'homologie de Hochschild $HH_*(A)$ de l'alg\`ebre~$A$.
Lorsque $N>2$, le couple $(C(A),b)$ n'est plus un complexe de cha\^\i nes.
On peut cependant d\'efinir des avatars ${}_p HH_*(A)$ de l'homologie de Hochschild
pour $p = 1, \ldots , N-1$ par
$${}_p HH_n(A) =
{\Ker \bigl( b^p : C_n(A) \to C_{n-p}(A) \bigr)
\over 
\Im \bigl( b^{N-p} : C_{n+N-p}(A) \to C_n(A) \bigr)
}. \eqno (0.2)$$
Notre premier r\'esultat (dont une version cohomologique
a \'et\'e annonc\'ee par Dubois-Violette [Dub])
\'enonce que, sous certaines conditions,
les groupes ainsi d\'efinis sont soit nuls,
soit isomorphes aux groupes de Hochschild usuels, 
ce qui nous donne une d\'efinition alternative et
exotique de l'homologie de Hochschild.

\medskip
\goodbreak
\noindent
{\sc Th\'eor\`eme~1.}---
{\it Consid\'erons un entier~$N\geq 2$, un scalaire $q$ et une alg\`ebre 
associative unif\`ere~$A$. 
Sous chacune des deux hypoth\`eses suivantes~:

(a) $q\neq 1$ est une racine primitive $N$-i\`eme de l'unit\'e,

(b) $N$ est un nombre premier, $q=1$ et $A$ est une $\ZZ/N\ZZ$-alg\`ebre, 

\noindent
nous avons, pour tout entier~$n\geq 0$ et tout $p=1,\ldots ,N-1$, 
$${}_p HH_n(A) \cong \left\{
\matrix{
HH_{2(n-p+1)/N}(A) & \quad \hbox{si} \;\; n+1 \equiv p \;\; \hbox{mod}\; N,\hfill\cr
\noalign{\medskip}
HH_{(2n+2-N)/N}(A) & \quad \hbox{si} \;\; n+1 \equiv 0 \;\; \hbox{mod}\; N,\hfill\cr
\noalign{\medskip}
0 & \quad \hbox{dans tous les autres cas.} \hfill\cr
}
\right.
$$
}
\medskip

Pour \'etablir ce th\'eor\`eme, 
nous nous pla\c cons dans le cadre g\'en\'eral des
$N$-complexes, c'est-\`a-dire des modules gradu\'es munis
d'un endomorphisme~$d$ de degr\'e~$-1$ tel que $d^N = 0$.
Ces objets ont \'et\'e abondamment utilis\'es dans les travaux r\'ecents 
de Dubois-Violette {\it et~al.\ } et de Kapranov 
sur le calcul diff\'erentiel quantique ({\it cf.}\ [Dub], [DK1], [DK2], [Kap]). 
Un embryon d'alg\`ebre homologique adapt\'ee aux $N$-complexes 
est apparu dans~[Kap].

Dans cet article, nous avons d\'evelopp\'e syst\'ematiquement
une telle alg\`ebre homologique en d\'efinissant la notion de $N$-r\'esolution
projective ou injective ainsi que  des groupes ${}_p\Tor$ et ${}_p\Ext$
qui jouissent de propri\'et\'es parall\`eles aux groupes $\Tor$ et $\Ext$
classiques. Le cadre dans lequel nous travaillons est
celui de l'alg\`ebre homologique relative d'Eilenberg et Moore.
Entre les groupes ${}_p\Tor$ et ${}_p\Ext$ et leurs pendants classiques,
nous \'etablissons des isomorphismes formellement similaires
\`a ceux du th\'eor\`eme~1.

Ce dernier est alors cons\'equence de l'isomorphisme
$${}_pHH_*(A,A) \cong {}_p\Tor_*^{A\ot A^{\op}}(A, A)
\eqno (0.3)$$
\'etabli au \S5 au moyen d'une $N$-r\'esolution projective
(relativement aux suites exactes de $A$-bimodules munies d'un
scindage $k$-lin\'eaire).
Pour d\'emontrer l'acyclicit\'e de cette r\'esolution, nous construisons
une homotopie \`a l'aide d'identit\'es remarquables dans l'alg\`ebre
des op\'erateurs aux $q$-diff\'erences engendr\'ee par deux ind\'etermin\'ees
$X$ et $Y$ soumises \`a la relation $YX-qXY = 1$. 
Ces identit\'es sont valides pr\'ecis\'ement
sous les m\^emes hypoth\`eses sur le scalaire~$q$ que celles du th\'eor\`eme~1.

Pour terminer, signalons que
des diff\'erentielles analogues \`a celles de~(0.1)
et des groupes d'homo\-logie du type (0.2) ont \'et\'e consid\'er\'es 
d\`es les ann\'ees~1940.
C'est ainsi que Mayer~[May] a introduit des groupes d'homologie singuli\`ere
des espaces topologiques \`a coefficients dans un corps fini $\ZZ/\ell$
avec une diff\'erentielle du type~(0.1) dans laquelle~$q=1$.
Spanier~[Spa] a d\'emontr\'e que les groupes d'homologie de Mayer se r\'eduisaient
aux groupes d'homologie traditionnels dans un th\'eor\`eme dont l'\'enonc\'e
est formellement le m\^eme que celui du th\'eor\`eme~1 plus haut.
Pour $q\neq 1$ il y a des r\'esultats
similaires de Sarkaria [Sar1],~[Sar2] concernant l'homologie simpliciale.

Nous remercions chaleureusement
Michel Dubois-Violette pour ses remarques et
ses encouragements.

\medskip\medskip\goodbreak
\noindent
{\sc Conventions et Notations.}
Dans toute la suite, on fixe un entier $N\geq 2$ et un anneau commutatif $k$ qu'on 
appellera l'anneau de base. Toutes les op\'erations d'alg\`ebre lin\'eaire
effectu\'ees le sont dans la cat\'egorie des $k$-modules.
Les $k$-alg\`ebres consid\'er\'ees dans ce travail sont associatives et unif\`eres.

On choisit un \'el\'ement $q$ de~$k$.
Pour tout $n\in \ZZ$, nous posons 
$$[n] = 
\left\{
\matrix{{q^n -1\over q-1} & \hbox{si} \; \; q\neq 1,\cr
\noalign{\medskip}
n & \hbox{si} \; \; q= 1.\cr
}\right.$$
Si $r\in \NN$ on convient que
$[0]! = 1$ et $[r]! = [1][2]\cdots [r]$ si $r>0$.
Enfin, pour $0\leq r,s \leq N-1$ 
on d\'efinit les $q$-coefficients binomiaux $(r,s)$ par
$$(r,s) = {[r+s]! \over [r]!\, [s]!} .$$

\vskip 25pt
\goodbreak

\vfill\eject

\noindent
{\sectionfont 1. G\'en\'eralit\'es 
sur les $N$-complexes et leur homologie}
\bigskip

\noindent
Nous nous pla\c cons dans une cat\'egorie ab\'elienne~$\cA$.

\medskip
\noindent
{\sc 1.1.\ D\'efinitions.}
Un {\it $N$-complexe} est un objet $\ZZ$-gradu\'e~$C$ de~$\cA$
muni d'un endomorphisme~$d$ de degr\'e~$-1$ tel que~$d^N = 0$. 
On notera $(C,d)$ la donn\'ee d'un $N$-complexe.
On dit que $(C,d)$ est {\it positif} (resp.\ {\it n\'egatif\ }) 
si $C_n = 0$ pour tout~$n<0$ (resp.\ pour tout~$n>0$).

Un {\it morphisme $f : (C,d) \to (C',d')$ de $N$-complexes} est un morphisme
$f : C \to C'$ d'objets $\ZZ$-gradu\'es de~$\cA$ qui est
de degr\'e~$0$ et tel que~$d' f = fd$.

Deux morphismes $f,g : (C,d) \to (C',d')$ de $N$-complexes
sont {\it homotopes} s'il existe un morphisme
$h : C \to C'$ de degr\'e~$N-1$ tel que
$$f-g = \sum_{i=0}^{N-1}\, d'{}^{N-1-i} h d^i.\eqno (1.1)$$
L'homotopie est une relation d'\'equivalence.

\medskip
\noindent
{\sc 1.2.\ Homologie.}
A tout $N$-complexe $(C,d)$ et tout $p=1,\ldots, N-1$ on peut associer
les {\it groupes d'homologie} ${}_p H_*(C)$ d\'efinis par
$${}_p H_n(C) =
{\Ker \bigl( d^p : C_n \to C_{n-p} \bigr)
\over 
\Im \bigl( d^{N-p} : C_{n+N-p} \to C_n \bigr)
}. \eqno (1.2)$$
Ces groupes sont fonctoriels, c'est-\`a-dire tout 
morphisme $f : (C,d) \to (C',d')$ de $N$-complexes induit
une application $f_* : {}_p H_*(C) \to {}_p H_*(C')$ telle
que $\id_*$ soit l'identit\'e et $(f\circ g)_* = f_* \circ g_*$
(lorsque $f$ et $g$ sont composables).

\medskip
\noindent
{\sc 1.3.\ Lemme.}---
{\it Deux morphismes homotopes $f, g : C\to C'$ de $N$-complexes induisent 
la m\^eme application sur les groupes d'homologie.}

\medskip
\Dem
Si $x\in C_n$ tel que $d^p x = 0$ repr\'esente un \'el\'ement
de~${}_pH_n(C)$, alors nous tirons de~(1.1) que
$$f(x) - g(x) = \sum_{i=0}^{p-1}\, d'{}^{N-1-i} h d^i(x)
= d'{}^{N-p}\Bigl( \sum_{i=0}^{p-1}\, d'{}^{p-1-i} h d^i(x) \Bigr) , $$
ce qui montre que $f$ et $g$ induisent la m\^eme application
${}_pH_n(C) \to {}_pH_n(C')$.
\hfill\cqfd

\medskip
\noindent
{\sc 1.4.\ Acyclicit\'e.}
On dit que le $N$-complexe $(C,d)$ est {\it acyclique} si
${}_p H_n(C) = 0$ pour tout $n\in \ZZ$ et tout $p=1, \ldots, N-1$.
Kapranov a montr\'e que pour que $(C,d)$ soit acyclique,
il suffit qu'il existe $p$ tel que 
${}_p H_n(C) = 0$ pour tout $n\in \ZZ$ (voir [Kap],~Proposition~1.5).

\medskip
\goodbreak
\noindent
{\sc 1.5.\ Les morphismes $i_*$ et $d_*$.}
Soit $(C,d)$ un $N$-complexe. Lorsque $N\geq 3$, on a deux
familles importantes de morphismes naturels
$$i_*:  {}_pH_i(C)\to  {}_{p+1}H_i(C) 
\quad\hbox{et}\quad
d_*:  {}_pH_i(C)\to  {}_{p-1}H_{i-1}(C).$$
La premi\`ere, de degr\'e~$0$, est induite par l'inclusion
$\Ker(d^p)\subset\Ker(d^{p+1})$.
La seconde, de degr\'e~$-1$, est induite par $d:\Ker(d^p)\to \Ker(d^{p-1})$.
On a le r\'esultat suivant dont on trouvera une d\'emonstration
dans [DK], Lemme~1.

\medskip
\noindent
{\sc 1.6.\ Lemme.}---
{\it
Pour tous les entiers strictement positifs $r$, $p$ tels que $0< p+r < N$, 
il~existe une longue suite exacte naturelle de la forme
$$\matrix{
{}_p H_*(C) & \hfl{i_*^r}{} & {}_{p+r} H_*(C) & \hfl{d_*^p}{} & {}_r H_*(C) \cr
\noalign{\smallskip}
\ufl{d_*^{N-p-r}}{} &&&& \vfl{i_*^{N-p-r}}{} \cr
\noalign{\smallskip}
{}_{N-r} H_*(C) & \rfl{i_*^p}{} & {}_{N-p-r} H_*(C) & 
\rfl{d_*^r}{} & {}_{N-p} H_*(C) \cr
}$$
}

\medskip
\noindent
{\sc 1.7.\ D\'efinition.}
Une {\it suite exacte courte 
$0\to C' \mapright{u} C \mapright{v}  C''\to 0$ 
de $N$-complexes}
est la donn\'ee de 
morphismes $C' \mapright{u}  C \mapright{v}  C''$ de $N$-complexes
tels que pour tout $n\in \ZZ$ la
suite $0\to C'_n \mapright{u}  C_n \mapright{v} C''_n\to 0$
est exacte.

\goodbreak
\medskip
\noindent
{\sc 1.8.\ Lemme.}
{\it Soit $0\to C' \mapright{u} C \mapright{v} C''\to 0$ une
suite exacte courte de $N$-complexes. 
Alors il existe une longue suite exacte naturelle
$$\matrix{
{}_p H_*(C') & \hfl{u_*}{} & {}_p H_*(C) & \hfl{v_*}{} & {}_p H_*(C'') \cr
\noalign{\smallskip}
\ufl{\del_2}{} &&&& \vfl{\del_1}{} \cr
\noalign{\smallskip}
{}_{N-p} H_*(C'') & \rfl{v_*}{} & {}_{N-p} H_*(C) & 
\rfl{u_*}{} & {}_{N-p} H_*(C') \cr
}$$
o\`u $\del_1$ est une application de degr\'e~$-p$ et 
$\del_2$ une application de degr\'e~$-(N-p)$.
}

\medskip
\Dem
Appliquons le lemme du serpent au diagramme de suites exactes
$$\matrix{
0 \to & C'_{n+N} & \mapright{u} & C_{n+N} & \mapright{v} &  C''_{n+N} & \to 0 \cr
\noalign{\smallskip}
& \vfl{d^{N-p}}{} && \vfl{d^{N-p}}{} && \vfl{d^{N-p}}{} & \cr
\noalign{\smallskip}
0\to & C'_{n+p}&  \mapright{u} & C_{n+p} & \mapright{v} & C''_{n+p} & \to 0 \cr
}$$
Nous obtenons pour tout~$n$ les deux suites exactes
$$0 \to \Ker(d^{N-p})\cap C'_n \; \mapright{u} \; \Ker(d^{N-p})\cap C_n
\; \mapright{v} \; \Ker(d^{N-p}) \cap C''_n$$
et 
$$C'_{n+p}/\Im(d^{N-p}) \; \mapright{u} \; C_{n+p}/\Im(d^{N-p}) 
\; \mapright{v} \; C''_{n+p}/\Im(d^{N-p}) \to 0, $$
que nous assemblons 
dans le diagramme commutatif de suites exactes
$$\matrix{
& C'_{n+p}/\Im(d^{N-p}) & \mapright{u} & C_{n+p}/\Im(d^{N-p}) & 
\mapright{v} & C''_{n+p}/\Im(d^{N-p}) & \to 0 \cr
\noalign{\smallskip}
& \vfl{d^{p}}{} && \vfl{d^{p}}{} && \vfl{d^{p}}{} & \cr
\noalign{\smallskip}
0 \to & \Ker(d^{N-p})\cap C'_n & \mapright{u} &  \Ker(d^{N-p})\cap C_n 
& \mapright{v} & \Ker(d^{N-p}) \cap C''_n & \cr
}$$
Une autre application du lemme du serpent nous donne la longue suite 
exacte d\'esir\'ee.
\hfill\cqfd

\goodbreak

\vskip 25pt

\noindent
{\sectionfont 2. R\'esolutions}
\bigskip

\noindent
On se place \`a nouveau dans une cat\'egorie ab\'elienne~$\cA$.

\medskip
\noindent
{\sc 2.1.\ Rappels d'Alg\`ebre Homologique Relative}
(voir [EM]).
Une {\it classe projective} dans $\cA$ est un couple $(\cP,\cS)$
o\`u $\cP$ est une classe d'objets de~$\cA$
et $\cS$~est une classe de suites $Y'\mapright{u} Y \mapright{v} Y''$ de~$\cA$ 
telles que~$vu=0$, v\'erifiant les trois axiomes suivants~:

(i) Un objet $P$ de $\cA$ est dans $\cP$ si et seulement si
pour toute suite $Y'\mapright{u} Y \mapright{v} Y''$ de~$\cS$,
la suite de groupes ab\'eliens
$$\Hom(P,Y') \to \Hom(P,Y) \to \Hom(P,Y'') $$
est exacte.

(ii) Une suite $Y'\mapright{u} Y \mapright{v} Y''$ de~$\cA$
est dans~$\cS$ si et seulement si pour tout objet $P$ de~$\cP$,
la suite
$$\Hom(P,Y') \to \Hom(P,Y) \to \Hom(P,Y'') $$
est exacte.

(iii) Pour tout morphisme $X\to Y$ de~$\cA$ 
il existe un objet $P$ de~$\cP$ et un morphisme $P\to X$ de~$\cA$
tels que $P\to X \to Y$ est dans~$\cS$.

Dans une classe projective $(\cP,\cS)$, la classe $\cP$ d'objets
et la classe $\cS$ de suites se d\'eterminent mutuellement.
Les objets de $\cP$ sont appel\'es $\cP$-projectifs.

\medskip
\noindent
{\sc 2.2.\ Exemples.}
(a) {\it (Classe projective absolue)} C'est la classe projective
$(\cP_{\abs},\cS_{\abs})$ pour laquelle $\cS_{\abs}$
est la classe de toutes les suites exactes de~$\cA$.

(b) {\it (Classe projective $k$-scind\'ee)}
Supposons que 
la cat\'egorie~$\cA$ soit $k$-lin\'eaire au-dessus de l'anneau de base~$k$.
Il existe alors une classe projective $(\cP_k,\cS_k)$
pour laquelle $\cS_k$ est constitu\'ee de toutes les suites
$Y'\mapright{u} Y \mapright{v} Y''$ telles que
le morphisme induit $Y' \to \Ker(v)$
est un \'epimorphisme $k$-scind\'e, c'est-\`a-dire poss\'edant un inverse \`a droite
qui est $k$-lin\'eaire.

Par exemple, si $\cA$ est la cat\'egorie des modules \`a gauche
sur une $k$-alg\`ebre~$A$, alors $\cP_k$~est form\'ee de tous
les $A$-modules de la forme $A\ot V$ o\`u $V$ est un $k$-module
quelconque.

\medskip\goodbreak
\noindent
{\sc 2.3.\ D\'efinitions.}
On fixe une classe projective $(\cP,\cS)$ de~$\cA$.
Un $N$-complexe $(C,d)$ est {\it $\cP$-exact} 
si le $N$-complexe induit $\Hom(P,C)$ est acyclique
pour tout objet $\cP$-projectif~$P$.

Soit $M$ un objet de la cat\'egorie~$\cA$. 
Une {\it $N$-r\'esolution $\cP$-projective} de~$M$
est la donn\'ee d'un $N$-complexe positif
$(P,d)$ de~$\cA$ et d'un morphisme $\eps : P_0\to M$ tels que

(i) l'objet $P_n$ est $\cP$-projectif pour tout $n\geq 0$,

(ii) la suite
$$\cdots \, \mapright{d} \, P_2 \, \mapright{d} \, P_1 \, \mapright{d} \, P_1
\, \mapright{\eps} \, M \to 0$$
est un $N$-complexe $\cP$-exact.

Nous noterons $(P,d,\eps :P_0\to M)$ la donn\'ee d'une $N$-r\'esolution
$\cP$-projective de~$M$.

\medskip
\noindent
{\sc 2.4.\ Existence de $N$-R\'esolutions $\cP$-Projectives.}
D\'emontrons que tout objet $M$ de~$\cA$ poss\`ede une $N$-r\'esolution
$\cP$-projective.
On sait que $M$ a une r\'esolution $\cP$-projective 
au sens de l'alg\`ebre homologique traditionnelle.
Choisissons-en une, soit $(Q, \delta) \mapright{\eps} M$.
Consid\'erons la suite d'applications
$$\cdots \, \mapright{\id} \, Q_2  \, \mapright{\delta} \, Q_1 
\, \mapright{\delta} \, Q_0 \, \mapright{\id}\, \cdots \, 
\mapright{\id} \, Q_0 \mapright{\eps} \, M \to 0
\eqno (2.1)$$
o\`u l'application identit\'e a \'et\'e ins\'er\'ee $N-2$~fois \`a la place de 
tous les objets~$Q_{2i}$~($i\geq 0$). 
On v\'erifie facilement que le $N$-complexe~(2.1) est une $N$-r\'esolution 
$\cP$-projective de~$M$.

\goodbreak\medskip
\noindent
{\sc 2.5.\ Contraction de $N$-Complexes.}
La construction du \S2.4 a une r\'eciproque.
Soit $(P,d)$ une $N$-r\'esolution $\cP$-projective d'un objet~$M$.
Etant donn\'e $p=1, \ldots, N-1$ nous lui associons la 
r\'esolution $\cP$-projective (traditionnelle) $\Delta_p P$
$$\cdots \mapright{d^{p}}\, P_{2N-1} \,\mapright{d^{N-p}}\, P_{N+p-1}
\, \mapright{d^{p}} \, P_{N-1} \, \mapright{d^{N-p}}\,  P_{p-1}
\, \mapright{\eps d^{p-1}}\,  M \to 0.$$
Notons qu'on a des morphismes de $N$-complexes
$$\Delta_{N-1}P\to \Delta_{N-2} P \to \cdots \to
\Delta_2 P \to \Delta_1 P  \eqno (2.2)$$ 
induits par la diff\'erentielle~$d$ comme suit~:
$$\matrix{
\cdots \mapright{d^{p+1}} & P_{2N-1} & \mapright{d^{N-p-1}}&  P_{N+p}& 
\mapright{d^{p+1}} & P_{N-1}&  \mapright{d^{N-p-1}}& P_{p}& 
\mapright{\eps d^{p}} & M &\to 0 \cr
\noalign{\smallskip}
& \vfl{\id}{} && \vfl{d}{} && \vfl{\id}{} && \vfl{d}{} && \vfl{\id}{} &\cr
\noalign{\smallskip}
\cdots \mapright{d^{p}} & P_{2N-1} & \mapright{d^{N-p}}&  P_{N+p-1} &
\mapright{d^{p}} &  P_{N-1} & \mapright{d^{N-p}} & P_{p-1}& 
\mapright{\eps d^{p-1}} & M &\to 0\cr
}$$
\goodbreak

Nous \'enon\c cons maintenant l'analogue pour les $N$-complexes
du lemme fondamental de l'alg\`ebre homologique.

\medskip
\noindent
{\sc 2.6.\ Lemme.}---
{\it Soit $u : M\to N$ un morphisme de~$\cA$. 
Si $(P,d,\eps :P_0\to M)$ est une $N$-r\'esolution $\cP$-projective de~$M$ et
$(Q,d,\eps :Q_0\to N)$ est une $N$-r\'esolution $\cP$-projective de~$N$, il existe
un morphisme de $N$-complexes $f: P\to Q$ tel que~$\eps f_0 = u\eps$.

Si $g: P\to Q$ est un autre morphisme de $N$-complexes
tel que~$\eps g_0 = u\eps$, alors $f$ et $g$ sont homotopes.
}

\medskip
\Dem
(a) Le $N$-complexe
$$\cdots \mapright{d}\, Q_1 \mapright{d}\, Q_0 \mapright{\eps}\, N \to 0
\eqno (2.3)$$
\'etant $\cP$-exact, l'application $\Hom(P_0,Q_0) \to \Hom(P_0,N)$
est surjective. Il existe donc un morphisme
$f_0 : P_0\to Q_0$ tel que $\eps f_0 = u\eps$.

Supposons construits des morphismes $f_i : P_i\to Q_i$ pour $i=0, \ldots, n$
tels que $d f_i = f_{i-1} d$  pour tout $1\leq i \leq n$.
Pour simplifier, posons $\eps = d$, $f_{-1} = u$ et 
$f_k = 0$ si $k < -1$.
On a ainsi $d^N = 0$ pour le $N$-complexe~(2.3) tout entier
et $d f_i = f_{i-1} d$  pour tout $i \leq n$.
L'hypoth\`ese de r\'ecurrence entra\^\i ne que 
$$d^{N-1} f_n d = f_{n-N+1} d^{N-1}d = f_{n-N+1} d^{N} = 0.$$
L'exactitude de la suite
$$\Hom(P_{n+1},Q_{n+1}) \to \Hom(P_{n+1},Q_n) \to \Hom(P_{n+1},Q_{n-N+1})$$
implique l'existence de 
$f_{n+1} : P_{n+1}\to Q_{n+1}$ tel que $d f_{n+1} = f_n d$.

(b) Si $g: P\to Q$ est un autre morphisme de $N$-complexes
tel que~$\eps g_0 = u\eps$, alors
$\eps(f_0 - g_0) = u\eps - u\eps = 0$.
L'exactitude de la suite
$$\Hom(P_0,Q_{N-1}) \to \Hom(P_0,Q_0) \to \Hom(P_0,N)$$ 
nous donne un morphisme $h_0 : P_0 \to Q_{N-1}$ tel que
$f_0 - g_0 = d^{N-1} h_0$.

Supposons construits des morphismes $h_k : P_k \to Q_{k+N-1}$ pour $k=0, \ldots, n$
tels que
$$f_k - g_k = \sum_{i=0}^{N-1}\, d^{N-i-1} h_{k-i} \, d^i$$
avec la convention $h_k = 0$ si $k<0$.
Alors
$$\eqalign{d(f_{n+1} - g_{n+1}) & = (f_n- g_n)d 
= \sum_{i=0}^{N-1}\, d^{N-i-1} h_{n-i} \, d^{i+1} \cr
& = \sum_{i=0}^{N-2}\, d^{N-i-1} h_{n-i} \, d^{i+1} 
= d\Bigl( \sum_{i=1}^{N-1}\, d^{N-i-1} h_{n+1-i} \, d^{i}\Bigr). \cr
}$$
Et donc $dF = 0$ avec
$$F = f_{n+1} - g_{n+1} - 
\Bigl( \sum_{i=1}^{N-1}\, d^{N-i-1} h_{n+1-i} \, d^{i}\Bigr).$$
La suite
$$\Hom(P_{n+1},Q_{n+N}) \to \Hom(P_{n+1},Q_{n+1}) \to \Hom(P_{n+1},Q_n)$$
\'etant exacte, il existe un \'el\'ement $h_{n+1}\in \Hom(P_{n+1},Q_{n+N})$ tel que
$F = d^{N-1}h_{n+1}$.
Par cons\'equent,
$$\eqalign{
f_{n+1} - g_{n+1} & =  d^{N-1}h_{n+1} + 
\sum_{i=1}^{N-1}\, d^{N-i-1} h_{n+1-i} \, d^{i}\cr
& = \sum_{i=0}^{N-1}\, d^{N-i-1} h_{n+1-i} \, d^{i}.\cr
}$$
On construit ainsi, de proche en proche, une homotopie entre $f$ et~$g$.
\hfill\cqfd

\vskip 25pt
\goodbreak

\noindent
{\sectionfont 3. Les groupes ${}_p\Tor_*^A(M,N)$}
\bigskip

\noindent
Soit $A$ une $k$-alg\`ebre et soit ${}_A\mod$ (resp.\ $\mod_A$)
la cat\'egorie ab\'elienne des $A$-modules \`a gauche (resp.\ \`a droite).
Dans chacune de ces cat\'egories, nous distinguons deux classes projectives,
\`a~savoir
la classe projective absolue $(\cP_{\abs},\cS_{\abs})$ et
la classe projective $k$-scind\'ee $(\cP_k,\cS_k)$.
Chaque classe projective fournit un bifoncteur gradu\'e~$\Tor_*^A(-,-)$.

Fixons l'une des classes projectives pr\'ec\'edentes $(\cP,\cS)$
et consid\'erons un $A$-module \`a droite~$M$ et
un $A$-module \`a gauche~$N$.

\medskip\goodbreak
\noindent
{\sc 3.1.\ Proposition.}---
{\it Soit $(P,d,\eps :P_0\to M)$ une $N$-r\'esolution $\cP$-projective de~$M$ et
$(Q,d,\eps :Q_0\to N)$ une $N$-r\'esolution $\cP$-projective de~$N$.
Pour tout $p=1,\ldots, N-1$ et tout $n\in \NN$, on a
$${}_p H_n(P\ot_A N, d\ot \id) = {}_p H_n(M\ot_A Q, \id\ot d)
\hskip 138pt
$$
$$\hskip 85pt
= \left\{
\matrix{
\Tor_{2(n-p+1)/N}^A(M,N) & \quad \hbox{si} \;\; n+1 \equiv p \;\; \hbox{mod}\; N,\hfill\cr
\noalign{\medskip}
\Tor_{(2n+2-N)/N}^A(M,N) & \quad \hbox{si} \;\; n+1 \equiv 0 \;\; \hbox{mod}\; N,\hfill\cr
\noalign{\medskip}
0 & \quad \hbox{sinon.} \hfill\cr
}
\right.
$$
}

\medskip
\noindent
{\sc 3.2.\ D\'efinition.}
La proposition pr\'ec\'edente montre que les groupes
${}_p H_*(P\ot_A N, d\ot \id)$ et ${}_p H_*(M\ot_A Q, \id\ot d)$
sont ind\'ependants des $N$-r\'esolutions $\cP$-projectives choisies.
Nous pouvons donc poser
$${}_p \Tor_n^A(M,N) = {}_p H_n(P\ot_A N) = {}_p H_n(M\ot_A Q)
\eqno (3.1)$$
o\`u $(P,d,\eps :P_0\to M)$ est une $N$-r\'esolution $\cP$-projective de~$M$
et $(Q,d,\eps :Q_0\to N)$ est une $N$-r\'esolution $\cP$-projective de~$N$.

\medskip
\noindent
{\sc 3.3.\ Corollaire.}---
{\it 
Les groupes ${}_p \Tor_n^A(M,N)$ sont li\'es aux groupes
$\Tor_n^A(M,N)$~par
$${}_p \Tor_n^A(M,N) = \left\{
\matrix{
\Tor_{2(n-p+1)/N}^A(M,N) & \quad \hbox{si} \;\; n+1 \equiv p \;\; \hbox{mod}\; N,\hfill\cr
\noalign{\medskip}
\Tor_{(2n+2-N)/N}^A(M,N) & \quad \hbox{si} \;\; n+1 \equiv 0 \;\; \hbox{mod}\; N,\hfill\cr
\noalign{\medskip}
0 & \quad \hbox{sinon.} \hfill\cr
}
\right.
$$
}

\medskip
\noindent
{\sc 3.4.\ D\'emonstration de la Proposition~3.1.}---
Soit  $p$  un entier compris entre 1 et $N-1$. 
Nous allons calculer ${}_p H_*(P\ot_A N, d\ot \id)$.
Les groupes ${}_p H_*(M\ot_A Q, \id\ot d)$ se calculent par la m\^eme m\'ethode.

Le complexe $\Delta_p P$ construit au \S2.5 \`a partir de la $N$-r\'esolution
$\cP$-projective est une r\'esolution $\cP$-projective de~$M$ au sens usuel.
Il en r\'esulte que
$\Tor_*^A(M,N)$ est l'homologie du complexe $\Delta_p P\ot_A N = \Delta_p (P\ot_A N)$~:
$$\cdots \mapright{(d\ot\id)^{p}} P_{2N-1}\ot_A N 
\mapright{(d\ot\id)^{N-p}} P_{N+p-1}\ot_A N
\mapright{(d\ot \id)^{p}} P_{N-1}\ot_A N \mapright{(d\ot\id)^{N-p}} P_{p-1}\ot_A N
\to 0.$$
Par cons\'equent, 
$$\Tor_{2n}^A(M,N) = {}_pH_{nN + p-1}(P\ot_A N) \et
\Tor_{2n+1}^A(M,N) = {}_{N-p}H_{(n+1)N -1}(P\ot_A N)$$
pour tout entier $n\geq 0$.
En particulier, le groupe ${}_{N-p}H_{(n+1)N -1}(P\ot_A N)$
ne d\'epend pas de~$p$, ce qui nous permet d'\'ecrire
$$\Tor_{2n+1}^A(M,N) = {}_{p}H_{(n+1)N -1}(P\ot_A N) .$$
On en d\'eduit aussit\^ot le calcul de ${}_p H_*(P\ot_A N, d\ot \id)$ lorsque 
$n+1\equiv 0$ ou $p$ modulo~$N$.
Si $N=2$ il n'y a plus rien \`a montrer.

Si $N>2$, il nous reste \`a \'etablir la nullit\'e de
${}_pH_{n}(P\ot_A N,d\ot\id)$ lorsque $n+1\not\equiv 0,p$ modulo~$N$.
Pour all\'eger les notations,  posons
$_pH_{n}(P\ot_A N,d\ot\id) = {}_pH_n$.

Du morphisme de r\'esolutions  $\Delta_{p+1} P\to\Delta_{p}P$ d\'efini au \S2.5,
nous d\'eduisons un 
quasi-isomorphisme $\Delta_{p+1} P\ot_ A N\to\Delta_{p} P\ot_ A N$.
En particulier, en degr\'e pair~$2k$, il induit l'isomorphisme
$$(d\ot\id)_* : {}_{p+1}H_{kN+p} \to {}_{p}H_{kN+p-1}$$
qui n'est autre que le morphisme $d_*$ du~\S1.5.
En degr\'e impair~$2k-1$, il induit l'isomorphisme
$$\id_* : {}_{N-p-1}H_{kN-1} \to {}_{N-p}H_{kN-1}$$
qui s'identifie avec le morphisme $i_*$ du~\S1.5.

Soit $r>0$ un entier tel que $p+r < N$.
Du lemme~1.6 nous tirons la suite exacte
$$\matrix{
{}_rH_{kN-1}\; \mapright{i_*^p} \; {}_{r+p}H_{kN-1} \; \mapright{d^r_*} 
& {}_{p}H_{kN-1-r} & \cr
\noalign{\medskip}
& \vfl{i^{N-r-p}_*}{} & \cr
\noalign{\medskip}
& {}_{N-r}H_{kN-1-r} & 
\mapright{d^p_*} \; {}_{N-r-p}H_{kN-1-r-p}.\cr
}$$
D'apr\`es ce qui pr\'ec\`ede, le premier morphisme
$i_*^p$ est un isomorphisme~; 
il en est de m\^eme du dernier morphisme $d^p_*$. 
On en d\'eduit que le
groupe $_{p}H_{kN-1-r}$ est nul.  En d'autres termes,
$_{p}H_{n}=0$ lorsque $n+1\equiv N-r$ modulo~$N$. 
Comme les conditions $r>0$ et $r+p<N$ sont \'equivalentes 
\`a la condition $p<N-r<N$, on a d\'emontr\'e que $_{p}H_{n}=0$ lorsque $n+1\equiv s$
avec~$p<s<N$.

Pour \'etablir que $_{p}H_{n}=0$ lorsque $n+1\equiv s$
avec~$0<s<p$, prenons l'entier $r$ tel que $0<r<p$ et $r+p<N$,
et consid\'erons la suite exacte \'egalement
tir\'ee du lemme~1.6
$$\matrix{
{}_{N-r} H_{(k+1)N-r-1} \; 
\mapright{d^{N-p}_*} \; {}_{p-r}H_{kN+p-r-1} \; \mapright{i^r_*} 
& {}_{p}H_{kN+p-r-1} & \cr
\noalign{\medskip}
& \vfl{d^{p-r}_*}{} & \cr
\noalign{\medskip}
& {}_{r}H_{kN-1} &
\mapright{i^{N-p}_*} \;  {}_{N-p+r}H_{kN-1}.\cr
}$$
Comme pr\'ec\'edemment, les morphismes extr\^emes sont des
isomorphismes. On en
d\'eduit que ${}_{p}H_{kN+p-r-1}=0$. Ainsi, on a prouv\'e que
$_{p}H_{n}=0$ lorsque $n+1\equiv p-r$ modulo~$N$.
On conclut en remarquant que $p-r$ parcourt tous les entiers $s$
tels que~$0<s<p$.
\hfill\cqfd

\medskip\goodbreak
\noindent
{\sc 3.5.\ Fonctorialit\'e.}
Les modules ${}_p \Tor_n^A(M,N)$ sont fonctoriels en $M$ et~$N$.
Cela r\'esulte du lemme~2.6 ou de l'identification donn\'ee
par le corollaire~3.3.

\medskip
\noindent
{\sc 3.6.\ Proposition.}
{\it (a) Soit $0\to N' \to N \to N''\to 0$ une
suite $\cP$-exacte courte de $N$-complexes de $A$-modules \`a gauche. 
Pour tout $A$-module \`a droite~$M$
il existe une longue suite exacte naturelle
$$\matrix{
{}_p \Tor_*(M,N') & \hfl{}{} & {}_p \Tor_*(M,N) & \hfl{}{} 
& {}_p \Tor_*(M,N'') \cr
\noalign{\smallskip}
\ufl{\del_2}{} &&&& \vfl{\del_1}{} \cr
\noalign{\smallskip}
{}_{N-p} \Tor_*(M,N'') & \rfl{}{} & {}_{N-p} \Tor_*(M,N) & \rfl{}{} 
& {}_{N-p} \Tor_*(M,N') \cr
}$$

(b) De m\^eme, si $0\to M' \to M \to M''\to 0$ est une
suite $\cP$-exacte courte de $N$-complexes de $A$-modules \`a droite
et si $N$ est un $A$-module \`a gauche, alors
il existe une longue suite exacte naturelle
$$\matrix{
{}_p \Tor_*(M',N) & \hfl{}{} & {}_p \Tor_*(M,N) & \hfl{}{} 
& {}_p \Tor_*(M'',N) \cr
\noalign{\smallskip}
\ufl{\del_2}{} &&&& \vfl{\del_1}{} \cr
\noalign{\smallskip}
{}_{N-p} \Tor_*(M'',N) & \rfl{}{} & {}_{N-p} \Tor_*(M,N) & \rfl{}{} 
& {}_{N-p} \Tor_*(M',N) \cr
}$$
Dans les cas (a) et (b), l'application $\del_1$ est de degr\'e~$-p$ et 
$\del_2$ de degr\'e~$-(N-p)$.
}

\medskip
\Dem
Nous nous contenterons d'\'etablir~(a). Le point (b) se d\'emontre de mani\`ere
similaire.

Soit $(P,d,\eps :P_0\to M)$ une $N$-r\'esolution $\cP$-projective de~$M$.
Consid\'erons la suite de $N$-complexes
$$0\to P\ot_A N' \to P\ot_A N \to P\ot_A N''\to 0 . \eqno (3.2)$$
La suite (3.2) est exacte \`a droite.
Si nous montrons que $P_n\ot_A N' \to P_n\ot_A N$ est injectif pour tout~$n\geq 0$, 
alors la suite (3.2) est exacte, ce qui permet de
lui appliquer le lemme~1.8 et de conclure.

Lorsque nous sommes dans le cas de la classe projective absolue,
le $N$-complexe est constitu\'e de $A$-modules qui sont facteurs directs
de modules libres. L'injectivit\'e de $P_n\ot_A N' \to P_n\ot_A N$
se ram\`ene au cas o\`u $P_n$ est un $A$-module libre, c'est-\`a-dire un cas
o\`u l'injectivit\'e est triviale.

Lorsque nous sommes dans la situation de la classe projective 
$k$-scind\'ee $(\cP_k,\cS_k)$, alors $P_n$ est de la forme $V\ot A$
avec la structure naturelle de $A$-module \`a droite.
L'exactitude de~(3.2) se r\'eduit alors \`a celle de la suite de $k$-modules
$$0\to V\ot N' \to V\ot N \to V\ot N''\to 0$$
pour tout $k$-module~$V$.
Cette derni\`ere est bien exacte. En effet, la $\cP_k$-exactitude de la
suite $0\to N' \to N \to N''\to 0$
implique qu'elle est exacte et scind\'ee dans la cat\'egorie des $k$-modules
(et r\'eciproquement),
propri\'et\'e qui est pr\'eserv\'ee par le foncteur~$V\ot -$.
\hfill\cqfd

\medskip\goodbreak
\noindent
{\sc 3.7.\ Caract\'erisation des Bifoncteurs ${}_p \Tor$.}
La famille de bifoncteurs gradu\'es ${}_p \Tor_*^A(-,-)$ ($0<p<N$)
est caract\'eris\'ee par les trois propri\'et\'es (i), (ii), (iii)
ci-dessous o\`u $M$ d\'esigne un $A$-module \`a droite et
$N$ un $A$-module \`a gauche.

(i) Pour tout $p$ tel que $0<p<N$, on a 
$${}_p \Tor_{p-1}^A(M,N) \cong M\ot _A N.$$

(ii) Si $P$ et $Q$ sont des modules $\cP$-projectifs 
et si $n\neq p-1$, on a
$${}_p \Tor_n^A(P,N) = 0 = {}_p \Tor_n^A(M,Q).$$

(iii) La proposition~3.6 est v\'erifi\'ee.

\noindent
Pour la d\'emonstration, on proc\`ede comme dans le cas classique.
\goodbreak

\vskip 25pt
\goodbreak

\noindent
{\sectionfont 4. Les groupes $\Ext_A^*(M,N)$}
\bigskip

\noindent
Soit $A$ une $k$-alg\`ebre et ${}_A\mod$
la cat\'egorie ab\'elienne des $A$-modules \`a gauche.
Nous nous restreignons volontairement aux
deux classes projectives $(\cP_{\abs},\cS_{\abs})$ et
$(\cP_k,\cS_k)$ consid\'er\'ees aux \S\S2--3.
A chacune de ces classes correspond la classe $\cI$ form\'ee
de tous les $A$-modules~$I$ tels que
la suite 
$$\Hom(Y'',I) \to \Hom(Y,I) \to \Hom(Y',I) $$
soit exacte 
pour toute suite $Y'\mapright{u} Y \mapright{v} Y''$ de~$\cS$. 
Un \'el\'ement de $\cI$ est appel\'e un module {\it $\cI$-injectif}. 

Ainsi, si $\cS= \cS_{\abs}$ est la classe de 
toutes les suites exactes de~${}_A\mod$, 
alors $\cI$ est constitu\'ee de ce que l'on appelle d'habitude
les $A$-modules injectifs.
Par contre, si $\cS = \cS_k$, alors $\cI$ est la classe des 
$A$-modules de la forme $\Hom_k(A,V)$ sur lesquels $A$
op\`ere par $(af)(a') = f(a'a)$
pour tout $f\in \Hom_k(A,V)$ et $a,a'\in A$.


\medskip
\noindent
{\sc 4.1.\ D\'efinition.}
Un $N$-complexe $(C,d)$ est {\it $\cI$-exact} 
si le $N$-complexe $\Hom(C,I)$ est acyclique
pour tout objet $\cI$-injectif~$I$.
 
Une {\it $N$-r\'esolution $\cI$-injective} d'un $A$-module~$M$
est la donn\'ee d'un $N$-complexe $\cI$-exact
$$0\to M \, \mapright{\eta}\, I_0 \, \mapright{d}\, I_{-1} \, \mapright{d}\,
I_{-2} \, \mapright{d}\, \cdots \eqno (4.1)$$
tel que le module $I_n$ est $\cI$-injectif pour tout $n\leq 0$.

Nous noterons $(I,d,\eta :M\to I_0)$ la donn\'ee d'une $N$-r\'esolution
$\cI$-injective de~$M$.

\medskip\goodbreak
\noindent
{\sc 4.2.\ Rapport avec les R\'esolutions Injectives Traditionnelles.}
L'existence de $N$-r\'esolutions $\cI$-injectives se d\'emontre comme
au \S2.4. On part d'une r\'esolution $\cI$-injective de~$M$
au sens de l'alg\`ebre homologique usuelle
$$0\to M \, \mapright{\eta}\, J_0 \, \mapright{d}\, J_{-1} \, \mapright{d}\,
J_{-2} \, \mapright{d}\, \cdots$$
On obtient une $N$-r\'esolution $\cI$-injective de~$M$ en ins\'erant
l'application identit\'e $N-2$~fois \`a la place de 
tous les modules~$J_{-2i}$~($i\geq 0$).

R\'eciproquement, si on part de la $N$-r\'esolution $\cI$-injective~(4.1),
on a, pour tout $p$ tel que $0<p<N$, 
la r\'esolution $\cI$-injective usuelle $\Delta^p J$
$$0\to M \, \mapright{d^{p-1}\eta}\, I_{-(p-1)} \, \mapright{d^{N-p}}\, 
I_{-(N-1)} \, \mapright{d^p}\,
I_{-(N+p-1)} \, \mapright{d^{N-p}}\, I_{-(2N-1)}
\, \mapright{d^p}\,\cdots \eqno (4.2)$$
Les r\'esolutions $\Delta^p I$ sont li\'ees par des morphismes de $N$-complexes
$$\Delta^{1}P\to \Delta^{2} P \to \cdots \to
\Delta^{N-2} P \to \Delta^{N-1} P  \eqno (4.3)$$
d\'efinis comme suit~:
$$\matrix{
0\to & \!\! M  \!\! & \!\! \mapright{d^{p-1}\eta} \!\! & \!\! I_{-(p-1)} \!\! 
& \!\! \mapright{d^{N-p}}\!\! & 
\!\! I_{-(N-1)} \!\! & \!\! \mapright{d^p} \!\! &
\!\! I_{-(N+p-1)} \!\! & \!\! \mapright{d^{N-p}} \!\! & \!\! I_{-(2N-1)}
\!\! & \!\! \mapright{d^p}\,\cdots\cr
\noalign{\smallskip}
& \vfl{\id}{} && \vfl{d}{} && \vfl{\id}{} && \vfl{d}{} && \vfl{\id}{} &\cr
\noalign{\smallskip}
0\to & M  & \mapright{d^{p}\eta} & I_{-p} & \mapright{d^{N-p-1}}& 
I_{-(N-1)} & \mapright{d^{p+1}} &
I_{-(N+p)} & \mapright{d^{N-p-1}} & I_{-(2N-1)}
& \mapright{d^{p+1}}\,\cdots\cr
}$$

Le lemme~2.6 a le pendant suivant dont nous laissons
la d\'emonstration au lecteur.

\medskip
\noindent
{\sc 4.3.\ Lemme.}---
{\it Soit $u : M\to N$ un morphisme de~$\cA$. 
Si $(I,d,\eta :M\to I_0)$ est une $N$-r\'esolution $\cI$-injective de~$M$ et
$(J,d,\eta :N\to J_0)$ est une $N$-r\'esolution $\cI$-injective de~$N$, il existe
un morphisme de $N$-complexes $f: I\to J$ tel que~$f_0\eta = \eta u$.

Si $g: I\to J$ est un autre morphisme de $N$-complexes
tel que~$g_0\eta = \eta u$, alors $f$ et $g$ sont homotopes.
}

\medskip
Nous fixons maintenant l'une des 
deux classes projectives $(\cP_{\abs},\cS_{\abs})$ et
$(\cP_k,\cS_k)$ de la cat\'egorie ab\'elienne~${}_A\mod$.
Il lui correspond une classe injective~$\cI$ et des
bifoncteurs $\Ext^*_A(-,-)$.

Consid\'erons des $A$-modules \`a gauche $M$ et~$N$.
Soit $(P,d,\eps :P_0\to M)$ une $N$-r\'esolution $\cP$-projective de~$M$ et
$(I,d,\eta :N\to I_0)$ une $N$-r\'esolution $\cI$-injective de~$N$.
Graduons $\Hom_A(M,I)$ et $\Hom_A(P,N)$ par
$$\Hom_A(M,I)_n = \Hom_A(M,I_n) \et \Hom_A(P,N)_n = \Hom_A(P_{-n},N)$$
et munissons-les des endomorphismes de degr\'e $-1$ donn\'es respectivement
par $f\mapsto d\circ f$ et $g\mapsto g\circ d$ pour tout 
$f\in \Hom_A(M,I)$ et $g\in \Hom_A(P,N)$.
On obtient ainsi des $N$-complexes n\'egatifs.

\medskip\goodbreak
\noindent
{\sc 4.4.\ Proposition.}---
{\it Sous les hypoth\`eses pr\'ec\'edentes, 
pour tout $p=1,\ldots, N-1$ et tout $n\in \NN$, on a
$${}_p H_{-n}(\Hom_A(M,I)) = {}_p H_{-n}(\Hom_A(P,N))
\hskip 138pt
$$
$$\hskip 85pt
= \left\{
\matrix{
\Ext^{2(n+1-(N-p))/N}_A(M,N) & \quad \hbox{si} \;\; n+1 \equiv N-p \;\; 
\hbox{mod}\; N,\hfill\cr
\noalign{\medskip}
\Ext^{(2n+2-N)/N}_A(M,N) & \quad \hbox{si} \;\; n+1 \equiv 0 \;\; 
\hbox{mod}\; N,\hfill\cr
\noalign{\medskip}
0 & \quad \hbox{sinon.} \hfill\cr
}
\right.
$$
}

\medskip
\Dem
Le complexe $\Delta_p P$ construit au \S2.5 est une r\'esolution 
$\cP$-projective de~$M$ au sens usuel. De m\^eme, 
$\Delta^p I$ est une r\'esolution $\cI$-injective de~$N$.
On peut donc calculer les groupes $\Ext^{*}_A(M,N)$
\`a partir des complexes $\Hom_A(M,\Delta^p I)$
et $\Hom_A(\Delta_p P,N)$. On obtient pour tout $n\in \NN$
$$\eqalign{
\Ext^{2n}_A(M,N) & = {}_{N-p} H_{-(nN+p-1)}(\Hom_A(M,I))
= {}_{p} H_{-(nN +N-p-1)}(\Hom_A(M,I)) \cr
& = {}_{N-p} H_{-(nN+p-1)}(\Hom_A(P,N))
= {}_{p} H_{-(nN +N-p-1)}(\Hom_A(P,N))\cr
}$$
et 
$$\Ext^{2n+1}_A(M,N) 
= {}_{p} H_{-[(n+1)N-1]}(\Hom_A(M,I)) 
= {}_{p} H_{-[(n+1)N-1]}(\Hom_A(P,N)).$$
On en d\'eduit, comme au~\S3.4, le calcul de ${}_p H_{-n}(\Hom_A(M,I))$ et de
${}_p H_{-n}(\Hom_A(P,N))$ lorsque $n+1\equiv 0$ ou $N-p$ modulo~$N$.

Quant aux autres cas, on les traite au moyen des morphismes
de r\'esolutions~(4.3) en proc\'edant comme au~\S3.4.
\hfill\cqfd

\medskip
\noindent
{\sc 4.5.\ D\'efinition.}
La proposition pr\'ec\'edente montre que les groupes
${}_p H_{*}(\Hom_A(M,I))$ et ${}_p H_{*}(\Hom_A(P,N))$
sont ind\'ependants des $N$-r\'esolutions choisies.
Nous pouvons donc poser
$${}_p \Ext^{n}_A(M,N) = {}_p H_{-n}(\Hom_A(M,I)) = {}_p H_{-n}(\Hom_A(P,N))$$
o\`u $(P,d,\eps :P_0\to M)$ est une $N$-r\'esolution $\cP$-projective de~$M$
et $(I,d,\eta :N\to I_0)$ est une $N$-r\'esolution $\cI$-injective de~$N$.

\medskip
\noindent
{\sc 4.6.\ Corollaire.}---
{\it 
Les groupes ${}_p \Ext^{n}_A(M,N)$ sont donn\'es~par
$${}_p \Ext^{n}_A(M,N) = \left\{
\matrix{
\Ext^{2(n+1-(N-p))/N}_A(M,N) & \quad \hbox{si} \;\; n+1 \equiv N-p \;\; 
\hbox{mod}\; N,\hfill\cr
\noalign{\medskip}
\Ext^{(2n+2-N)/N}_A(M,N) & \quad \hbox{si} \;\; n+1 \equiv 0 \;\; 
\hbox{mod}\; N,\hfill\cr
\noalign{\medskip}
0 & \quad \hbox{sinon.} \hfill\cr
}
\right.
$$
}

\medskip\goodbreak
\noindent
{\sc 4.7.\ Caract\'erisation des Bifoncteurs ${}_p \Ext$.}
Les modules ${}_p \Ext^{n}_A(M,N)$ sont fonctoriels en 
les $A$-modules \`a gauche $M$ et~$N$.
Ils sont caract\'eris\'es par les propri\'et\'es (i)--(iv)
ci-dessous.

(i) Pour tout $p$ tel que $0<p<N$, on a
$${}_p \Ext^{N-p-1}_A(M,N) \cong \Hom_A(M,N).$$

(ii) Pour tous $n$ et $p$ tels que $n\neq N-p-1$, tout module $\cP$-projectif~$P$
et tout module $\cI$-injectif~$I$, on a
$${}_p \Ext^{n}_A(P,N) = 0 = {}_p \Ext^{n}_A(M,I).$$

(iii) Soit $0\to N' \to N \to N''\to 0$ une
suite $\cP$-exacte courte de $N$-complexes de $A$-modules. 
Pour tout $A$-module~$M$
il existe une longue suite exacte naturelle
$$\matrix{
{}_p \Ext_A^*(M,N') & \hfl{}{} & {}_p \Ext_A^*(M,N) & \hfl{}{} 
& {}_p \Ext_A^*(M,N'') \cr
\noalign{\smallskip}
\ufl{\del_2}{} &&&& \vfl{\del_1}{} \cr
\noalign{\smallskip}
{}_{N-p} \Ext_A^*(M,N'') & \rfl{}{} & {}_{N-p} \Ext_A^*(M,N) & \rfl{}{} 
& {}_{N-p} \Ext_A^*(M,N') \cr
}$$

(iv) De m\^eme, si $0\to M' \to M \to M''\to 0$ est une
suite $\cI$-exacte courte de $N$-complexes de $A$-modules
et si $N$ est un $A$-module quelconque, alors
il existe une longue suite exacte naturelle
$$\matrix{
{}_p \Ext_A^*(M'',N) & \hfl{}{} & {}_p \Ext_A^*(M,N) & \hfl{}{} 
& {}_p \Ext_A^*(M',N) \cr
\noalign{\smallskip}
\ufl{\del_2}{} &&&& \vfl{\del_1}{} \cr
\noalign{\smallskip}
{}_{N-p} \Ext_A^*(M',N) & \rfl{}{} & {}_{N-p} \Ext_A^*(M,N) & \rfl{}{} 
& {}_{N-p} \Ext_A^*(M'',N) \cr
}$$
Les applications $\del_1$ et $\del_2$ ci-dessus sont de degr\'es
respectifs $p$ et~$N-p$.

Dans la d\'emonstration de l'exactitude des longues suites
pr\'ec\'edentes pour le cas o\`u $\cS = \cS_k$,
on utilise les isomorphismes d'adjonction 
$$\Hom_A(A\ot_k V,M) \cong \Hom_k(V,M) \et
\Hom_A(M,\Hom_k(A,V)) \cong \Hom_k(M,V)$$
o\`u $M$ est un $A$-module \`a gauche et $V$ un $k$-module.

%

\vskip 25pt
\goodbreak

\noindent
{\sectionfont 5. D\'emonstration du th\'eor\`eme~1}
\bigskip

\noindent
Dans ce paragraphe nous aurons besoin des deux hypoth\`eses
(H$_0$) et (H$_1$) suivantes, d\'ej\`a mentionn\'ees dans l'introduction.
Elles portent sur l'entier $N\geq 2$, l'anneau commutatif~$k$ et 
l'\'el\'ement~$q$ de~$k$.

\smallskip
\noindent
(H$_0$) On a $[N] = 0$ dans $k$.

\smallskip
\noindent
(H$_1$) L'hypoth\`ese H$_0$ est v\'erifi\'ee et,
de plus, $[p]$ est inversible dans~$k$ pour tout~$0<p<N$.
\smallskip

L'hypoth\`ese (H$_0$) est \'equivalente \`a~:
soit $q$ est une racine $N$-i\`eme de l'unit\'e diff\'erente de~$1$,
soit $q=1$ et $N=0$ dans~$k$.

L'hypoth\`ese (H$_1$) est \'equivalente \`a~:
soit $q$ est une racine primitive $N$-i\`eme de l'unit\'e,
soit $q=1$, l'entier $N$ est premier et $k$ est une alg\`ebre
sur le corps~$\ZZ/N$.

\medskip\goodbreak
\noindent
{\sc 5.1.\ Modules simpliciaux.}
Soit $(C, d_i, s_i)$ un $k$-module simplicial muni d'applications
face $d_i$ et d\'eg\'en\'erescence~$s_i$. On a notamment
les formules de commutation simpliciales
$$d_id_j = d_{j-1}d_i \quad\hbox{si}\; i<j  \et
d_is_j = \left\{\matrix{
s_{j-1}d_i \hfill & \hbox{si}\; i<j, \hfill\cr
\noalign{\smallskip}
\id \hfill & \hbox{si}\; i=j, j+1, \hfill\cr
\noalign{\smallskip}
s_jd_{i-1} \hfill & \hbox{si}\; i>j+1. \hfill\cr
}\right.
$$
Pour tout scalaire $q\in k$, on d\'efinit
$d : C_n\to C_{n-1}$ par
$$d = \sum_{i=0}^n q^i\, d_i. \eqno (5.1)$$
Kapranov ([Kap], Proposition~0.2) et
Dubois-Violette ([Dub], Lemme~3)
ont observ\'e que l'on a $d^N = 0$ sous l'hypo\-th\`ese~(H$_0$).

\medskip
\noindent
{\sc 5.2.\ Proposition.}---
{\it Avec les notations pr\'ec\'edentes, d\'efinissons
$\delta : C_n\to C_{n-1}$ par
$$\delta = \sum_{i=0}^{n-1} q^i\, d_i. \eqno (5.2)$$
On a $\delta^N = 0$ sous l'hypoth\`ese~(H$_0$)~;
si, de plus, l'hypoth\`ese~(H$_1$) est v\'erifi\'ee, alors
$(C,\delta)$ est un $N$-complexe acyclique.
}
\goodbreak\medskip

Nous allons nous placer dans la situation plus g\'en\'erale suivante.
Supposons que $\delta : C_n\to C_{n-1}$ soit de la forme
$$\delta = \sum_{i=0}^{n-1} a_{n-1-i}\, q^i\, d_i $$
o\`u $a_0, a_1, \ldots$ sont des scalaires fix\'es une fois pour toutes.
\goodbreak

\medskip
\noindent
{\sc 5.3.\ Lemme.}---
{\it Si l'on pose $(N,j,k) = \sum_{s=0}^{N-k}\, q^sa_{n-k-j-s}$, on a
$$\delta^N = \sum_{0\leq i_1\leq i_2\leq \cdots \leq i_N\leq n-N}\, 
q^{i_1 + i_2 + \cdots + i_N} \prod_{k=1}^N \, (N,i_k,k) \;
d_{i_N}d_{i_{N-1}} \ldots d_{i_2}d_{i_1}.$$
}

\medskip
\Dem
Nous allons proc\' eder par r\' ecurrence sur l'entier~$N$.
La formule est  v\' erifi\' ee lorsque~$N=1$.
Supposons qu'elle le soit pour un entier $N\geq 1$.
Nous avons
$$\eqalign{
\delta^{N+1} & =\delta\delta^{N}
=\Bigl( \sum_{j=0}^{n-N-1} \, a_{n-N-1-j}\, q^jd_{j}\Bigr) \delta^N \cr
& = \sum_{0\le i_1\le\cdots\le i_N\le n-N\atop 0\le j\le n-N-1}
\hskip-15pt
q^{i_1+\cdots+i_N+j}\, a_{n-N-1-j}\,
\prod_{k=1}^N \, (N,i_k,k) \; d_jd_{i_N}\ldots d_{i_1} \cr
& = \sum_{r=0}^N\, S_r\cr
}$$
o\`u $S_r$ est la somme des termes pr\'ec\'edents
pour lesquels on a l'encadre\-ment $i_r\leq j< i_{r+1}$ si $0\leq r \leq N-1$
(convenons que $i_0=0$) et $i_N\leq j$ si~$r=N$.
Des relations de commutation simpliciales on tire pour $r<N$
$$S_r =
\sum \, q^{i_1+\cdots+i_N+j}\, a_{n-N-1-j}\, \prod_{k=1}^N\, (N,i_k,k)\;
d_{i_N-1}\ldots d_{i_{r+1}-1}d_jd_{i_r}\ldots d_{i_1}.
$$
En renommant les indices, on obtient
$$S_r = 
\sum q^{i_1+\cdots+i_{N+1}}q^{N-r}\, a_{n-N-1-i_{r+1}}\,
\prod_{k=1}^{r}\, (N,i_k,k)\, \prod_{k=r+1}^{N}\, (N,i_{k+1}+1,k)\,
d_{i_{N+1}}d_{i_N}\ldots d_{i_1}
$$
o\`u la somme est prise sur $0\le i_1\le\cdots\le i_{N+1}\le n-N-1$.
Posons
$$\Phi_r = q^{N-r}\, a_{n-N-1-i_{r+1}}\,
\prod_{k=1}^{r}\, (N,i_k,k)\, \prod_{k=r+1}^{N}\, (N,i_ {k+1}+1,k)$$
si $0\leq r\leq N-1$~; sinon,
$$\Phi_{N} = (N+1,i_{N+1},N+1)\, \prod_{k=1}^{N}\, (N,i_{k},k).$$
\goodbreak\noindent
En sommant sur $r=0,\ldots, N$, on a 
$$\delta^{N+1} =
\sum_{0\le i_1\le\cdots\le i_{N+1}\le n-N-1}\hskip-20pt
q^{i_1+\cdots+i_{N+1}} 
\, \Bigl(\sum_{r=0}^{N}\Phi_r\Bigr) \,
d_{i_{N+1}}\ldots d_{i_1}.
$$
Il ne reste plus qu'\` a montrer que
$$\sum_{r=0}^{N}\, \Phi_r = \prod_{k=1}^{N+1}\, (N+1,i_k,k).$$
Remarquons que, pour tous les entiers $i$, $k$,
$$(N+1,i,k)-(N,i,k)=q^{N+1-k}a_{n-1-N-i} \et
(N,i+1,k) = (N+1,i,k+1).$$
Par cons\' equent,
$$\eqalign{
\Phi_r & = \bigl[ (N+1,i_{r+1},r+1)-(N,i_{r+1},r+1)\bigr]\,
\prod_{k=1}^{r}\, (N,i_{k},k)\,
\prod_{k=r+1}^{N}\, (N+1,i_{k+1},k+1) \cr
& =\prod_{k=1}^{r}\, (N,i_{k},k)\,
   \prod_{k=r+1}^{N+1}\, (N+1,i_{k},k) -
    \prod_{k=1}^{r+1}\, (N,i_{k},k)\,
    \prod_{k=r+2}^{N+1}\, (N+1,i_{k},k) .\cr}
$$
Les termes de la somme $\sum_{r=0}^{N}\Phi_r$ s'\'eliminent deux \`a deux
sauf $\prod_{k=1}^{N+1}\, (N+1,i_k,k)$.
\hfill\cqfd

\medskip\goodbreak
\noindent
{\sc 5.4.\ Corollaire.}---
{\it D\'efinissons $\delta : C_n\to C_{n-1}$ par
$$\delta = \sum_{i=0}^{n-2} q^i\, d_i  + [\ell]\,  q^{n-1} d_{n-1} \eqno (5.3)$$
o\`u $\ell$ est un entier quelconque. Alors $\delta^N=0$
sous l'hypoth\`ese~(H$_0$).
}

\medskip
\Dem
Comme $[\ell]=[m]$ lorsque $\ell\equiv m$ modulo~$N$,
nous pouvons nous limiter au cas o\`u $0\leq \ell\leq N-1$.

Supposons dans un premier temps que $\ell=1$.
Nous appliquons le lemme~5.3 au cas o\`u 
$a_i=1$ pour tout~$i\geq 0$.
Le produit $\prod_{k=1}^N\, (N,i_k,k)$
est nul, et donc $\delta^N=0$, car
$$(N,i_1,1) = \sum_{s=0}^{N-1}\, q^{s}\, a_{n-1-i_1-s}
= \sum_{s=0}^{N-1}\, q^{s} = [N] = 0 .$$

Si maintenant $\ell=0$, alors $\delta = \sum_{i=0}^{n-2}\, q^i d_i$. 
En rempla\c cant $n$ par~$n-1$, on se ram\`ene
au cas $\ell=1$.

Si $\ell>1$ on applique le lemme~5.3 avec $a_0 = [\ell]$ et
$a_i=1$ pour tout~$i\geq 1$.
Montrons \`a nouveau que le produit $\prod_{k=1}^{N}\, (N,i_k,k)$ 
est nul en \'etablissant la nullit\'e d'un de ses facteurs.

Lorsque $i_1 < n-N$, on a $a_{n-N-i_1+r} = 1$ pour $r\ge0$. 
Alors le facteur $(N,i_1,1)$ est nul
en vertu du m\^eme calcul que dans le cas~$\ell=1$.

Lorsque $i_1=n-N$, il en est {\it a fortiori} de m\^eme de~$i_{\ell}$. 
Maintenant, c'est le facteur $(N,i_\ell,\ell)$ qui est nul. En effet,
$$\eqalign{
(N,i_\ell,\ell) &= \sum_{s=0}^{N-\ell}\, q^{s}\, a_{n-\ell -i_{\ell}-s} 
= \sum_{s=0}^{N-\ell}\, q^{s}\, a_{N-\ell -s} \cr
& = q^{N-\ell}\, [\ell] + \sum_{s=0}^{N-\ell-1}\, q^{s} 
=q^{N-\ell}\, [\ell] + [N-\ell] \cr
& = q^{N-\ell} \, [\ell] + [-\ell]=0.\cr}
$$
En conclusion, $\delta^N=0$.
\hfill\cqfd
\medskip
\goodbreak

Avant de passer \`a la d\'emonstration de la Proposition~5.2,
d\'emontrons deux relations dans
la $k$-alg\`ebre  $\cD_q$ engendr\'ee par les variables
$X$ et $Y$ et la relation $YX - qXY=1$.

\medskip
\noindent
{\sc 5.5.\ Lemme.}--- 
{\it Sous l'hypoth\`ese (H$_1$) on a les relations suivantes dans $\cD_q$~:
$$\sum_{k=0}^{N-1}\, X^{N-k-1}Y^{N-1}X^k = [N-1]!,$$
$$\sum_{k=0}^{N-1}\, Y^{N-k-1}X^{N-1}Y^k = (-1)^{N-1} q^{-N(N-1)/2}[N-1]!. $$
}
\medskip

\Dem
Commen\c cons par \'etablir la premi\`ere relation.
Remarquons tout d'abord que dans $\cD_q$ on a la formule\footnote{$^1$}{
\notefont
L'alg\`ebre $\scriptstyle \cD_q$ s'obtient aussi comme extension d'Ore de
l'alg\`ebre de polyn\^omes~$\scriptstyle k[X]$ \`a l'aide des op\'erateurs 
$\scriptstyle \del_q$ et $\scriptstyle \tau_q$ d\'efinis plus bas. 
La formule~(5.4) est un cas particulier
d'identit\'es du m\^eme type valides dans toute extension d'Ore.}
$$ Y^{\ell}X^{k} =
\sum_{r=0}^{\ell} \, q^{(\ell-r)(k-r)} \, (\ell-r,r)\,
[k][k-1]\cdots[k-r+1]X^{k-r}Y^{\ell-r} \eqno (5.4)$$
pour tout $k,\ell\geq 0$.
Elle peut se d\'emontrer par r\'ecurrence sur~$\ell$.
Pour $\ell=0$, il n'y a rien \`a \'etablir.
Pour $\ell=1$ la formule (5.4) a la forme
$$YX^k = q^k X^kY + [k] X^{k-1}, \eqno (5.5)$$
formule qu'on \'etablit par r\'ecurrence sur~$k$.
Pour $\ell>1$, on calcule $Y^{\ell}X^{k} = Y(Y^{\ell-1}X^{k})$
\`a partir de la formule (5.4) pour $Y^{\ell-1}X^{k}$ et de~(5.5).

En utilisant (5.4) on obtient
$$\sum_{k=0}^{N-1}\, X^{N-k-1}Y^{N-1}X^k  \hskip 270pt$$
$$\eqalign{
&= \sum_{k=0}^{N-1}\, 
\sum_{r=0}^{N-1} \, q^{(N-1-r)(k-r)}\, (N-1-r,r)\, 
[k][k-1]\cdots [k-r+1] \, X^{N-1-r}Y^{N-1-r}\cr
&= \sum_{r=0}^{N-1}\,  (N-1-r,r)\, a(r) \, X^{N-1-r}Y^{N-1-r}\cr
}$$ 
o\`u
$$a(r) = \sum_{k=r}^{N-1} \, q^{(N-1-r)(k-r)}\, [k][k-1]\cdots [k-r+1].$$
Nous affirmons que $a(r) = 0$ pour $0\leq r< N-1$.
Il en r\'esulte que
$$\sum_{k=0}^{N-1}\, X^{N-k-1}Y^{N-1}X^k = (0,N-1)\, a(N-1) = [N-1]!,$$
ce qui d\'emontre la premi\`ere identit\'e du lemme~5.5.

\goodbreak
Pour \'etablir la nullit\'e de $a(r)$ lorsque $0\leq r< N-1$,
introduisons sur l'alg\`ebre de polyn\^omes $k[X]$
l'op\'erateur aux $q$-diff\'erences $\del_q$ d\'efini par 
$$\del_q(X^k) = [k]\,  X^{k-1}$$
pour tout $k\geq 0$. Alors
$$\eqalign{
a(r) & = 
\sum_{k=r}^{N-1}\, [k][k-1]\cdots[k-r+1]\, X^{k-r}{}_{\textstyle |_{X=q^{N-1-r}}}\cr
&=\del_q^r \left(\sum_{k=0}^{N-1}\, X^{k} \right)\!_{\textstyle |_{X=q^{N-1-r}}} \cr
& = \del_q^r \left({X^N-1\over X-1}\right)\!_{\textstyle |_{X=q^{N-1-r}}}
.\cr
}$$
Admettons pour l'instant que
$$\del_q^r \left({X^N-1\over X-1}\right) =
(-1)^r{[r]!(X^N-1)\over (X-1)(qX-1)\cdots(q^rX-1)}. \eqno (5.6)$$

Si $q=1$, alors $N$ est un nombre premier sous l'hypoth\`ese~(H$_1$)
et donc
$${(X^N-1)\over (X-1)(qX-1)\cdots(q^rX-1)}
= {(X-1)^N\over (X-1)^{r+1}} = (X-1)^{N-r-1},$$
polyn\^ome qui admet $1$ comme racine lorsque $r<N-1$.
Par cons\'equent,
$$a(r) = (-1)^r r!\, (X-1)^{N-r-1}\!_{\textstyle |_{X=1}}
=0.$$
\goodbreak

Si $q\neq 1$, alors $q$ est une racine primitive $N$-i\`eme de 1,
ce qui nous permet d'\'ecrire $X^N -1$ comme produits des 
polyn\^omes $(X-q^{-i})$ o\`u $i=0, \ldots ,N-1$.
Il en r\'esulte que 
$${(X^N-1)\over (X-1)(qX-1)\cdots(q^rX-1)} = 
q^{-r(r+1)/2}\, \prod_{i=r+1}^{N-1}\, (X-q^{-i})$$
s'annule pour $X=q^{N-1-r} = q^{-(r+1)}$. 
La nullit\'e de $a(r)$ pour~$r<N-1$ en d\'ecoule.

Il reste \`a d\'emontrer~(5.6). Posons
$$P_r = \del_q^r \left({X^N-1\over X-1}\right).$$
On va \'etablir~(5.6) par r\'ecurrence sur~$r$.
C'est vrai pour $r=0$.
On sait que l'op\'erateur $\del_q$ est une $\tau_q$-d\'erivation de $k[X]$
o\`u $\tau_q$ est l'automorphisme d'alg\`ebre qui envoie $X$ sur $qX$, 
c'est-\`a-dire que l'on a 
$$\del_q(P_1P_2\ldots P_n) = 
\sum_{i=1}^n\, \tau_q(P_1\ldots  P_{i-1})\, \del_q(P_i) \, P_{i+1}\ldots P_n
\eqno (5.7)$$
pour tous $P_1, P_2, \ldots, P_n \in k[X]$.
Supposons la formule~(5.6) d\'emontr\'ee jusqu'au cran~$r$.
Appliquons l'op\'erateur $\del_q$ aux deux membres de l'identit\'e
$$(-1)^r [r]!(X^N-1) = (X-1)(qX-1)(q^rX-1)P_r.$$
A gauche, on obtient $(-1)^r [r]! [N] X^{N-1}$
qui est nul sous l'hypoth\`ese~(H$_1$).
A droite, la formule de Leibniz~(5.7) donne
$$(qX-1)\ldots (q^rX-1) \Bigl( (1+ q + \cdots + q^r)\, P_r + (q^{r+1}X-1)P_{r+1}
\Bigr).$$
Par cons\'equent,
$$P_{r+1} = - {[r+1]P_r\over (q^{r+1}X-1)} = 
(-1)^{r+1}{[r+1]!(X^N-1)\over (X-1)(qX-1)\cdots(q^rX-1)(q^{r+1}X-1)},$$
ce qu'il fallait d\'emontrer.
 
La seconde formule du lemme~5.5 se ram\`ene \`a la premi\`ere comme suit.
Consid\'erons l'alg\`ebre~$\cD'_q$ engendr\'ee  par  les
variables $\tilde X$ et $\tilde Y$ et 
la relation $\tilde Y\tilde X-q^{-1} \tilde X \tilde Y =1$.
Elle est obtenue de~$\cD_q$ en rempla\c cant $q$ par~$q^{-1}$.
Dans~$\cD'_q$ la premi\`ere formule du lemme~5.5 s'\'ecrit
$$\sum_{k=0}^{N-1}\, \tilde X^{N-k-1}\tilde Y^{N-1}\tilde X^k
= q^{-(N-1)(N-2)/2}[N-1]! . \eqno (5.8)$$
Il ne reste plus qu'\`a lui appliquer le morphisme d'alg\`ebres
$\alpha:\cD'_q\to\cD_q$ d\'efini par 
$$\alpha(\tilde X)=-qY \et \alpha(\tilde Y)=X,$$ 
ce qui nous donne la seconde formule du lemme. 
\hfill\cqfd

\medskip\goodbreak
\noindent
{\sc 5.6.\ D\'emonstration de la Proposition~5.2.}---
L'assertion $\delta^N=0$ r\'esulte du corollaire~5.4 avec $\ell=1$.

Etablissons maintenant l'acyclicit\'e de $(C,\delta)$.
Consid\'erons l'application $\sigma:C_{n}\to C_{n+1}$
d\'efinie par $\sigma = q^{-n}s_n$. On a 
$$\delta\sigma-q^{-1}\sigma\delta=\id.$$
En effet, les relations de commutation simpliciales donnent
$$\eqalign{
\delta\sigma 
& =\sum_{i=0}^n\, q^{i-n}d_is_n 
=\sum_{i=0}^{n-1}\, q^{i-n}d_is_n+ d_ns_n \cr
& = \sum_{i=0}^{n-1}\, q^{i-n}s_{n-1}d_i+ \id 
=q^{-n}s_{n-1} \Bigl(\sum_{i=0}^{n-1}\, q^{i}d_i\Bigr) + \id \cr
& =q^{-1}\sigma\delta+\id.\cr
}$$
Sous l'hypoth\`ese~(H$_1$) on peut appliquer la premi\`ere formule
du lemme~5.5
avec $X = \delta$ et $Y = -q^{-1}\sigma$,
ce qui donne
$$(-1)^{N-1}q^{-(N-1)}\sum_{k=0}^{N-1}\, \delta^{N-k-1}\sigma^{N-1}\delta^k 
= [N-1]!\, \id.$$
Puisque $[N-1]!$ est inversible sous nos hypoth\`eses, l'application
$\sigma^{N-1}$ convenablement normalis\'ee est une homotopie entre l'identit\'e et 
l'application nulle.
\hfill\cqfd

\medskip
\noindent
{\sc 5.7.\ Remarque.}---
Le $N$-complexe $C$ muni de la diff\'erentielle~(5.3)
est encore acyclique lorsqu'il existe une application lin\'eaire
$s : C_n\to C_{n+1}$ telle que
$$d_0 s = \id \quad\hbox{et}\quad d_is = s d_{i-1}$$
pour tout $i>0$.
En effet, dans ce cas on a 
$$\eqalign{
\delta s&=\sum_{i=0}^{n-1}\, q^i\, d_is+[\ell]\, q^n\,  d_ns
= d_0s+\sum_{i=1}^{n-1}\, q^i\, d_is + [\ell]\, q^n\, d_ns \cr
& = \id+\sum_{i=1}^{n-1}\, q^i\, sd_{i-1} + [\ell]\, q^n\, sd_{n-1}\cr
&=\id + \sum_{i=0}^{n-2}\, q^{i+1}\, sd_{i} + 
[\ell]\, q^n\, sd_{n-1} = \id + q\, s\delta.\cr
}$$
La relation $\delta s - q s \delta = \id$ 
nous permet d'appliquer la seconde formule du lemme~5.5
avec $X = s$ et $Y = \delta$, et d'obtenir ainsi l'identit\'e
$$\sum_{k=0}^{N-1}\, \delta^{N-k-1}\, s^{N-1}\, \delta^k 
= (-1)^{N-1} q^{-N(N-1)/2}[N-1]! \, \id.$$
On conclut comme au~\S5.6.

Ceci s'applique, en particulier, au cas $\ell = -1$, o\`u l'on retrouve
la diff\'erentielle $d_q$ consid\'er\'ee dans~[Dub].
De m\^eme, le
lemme~2 de [Dub] qui d\'eduit l'acyclicit\'e de $(C,\delta)$
de l'existence d'une application $s$ v\'erifiant
$s\delta  - q  \delta s = \id$, peut se d\'emontrer
\`a l'aide de la premi\`ere identit\'e du lemme~5.5.

\medskip\goodbreak
\noindent
{\sc 5.8.\ Application.}
Soit $A$ une alg\`ebre associative unif\`ere
et $(C(A), d_i, s_i)$ le module simplicial
de Hochschild de $A$. Rappelons que $C_n(A) = A^{\ot (n+1)}$
et que les applications $d_i : C_n(A) \to C_{n-1}(A)$
et $s_i : C_n(A) \to C_{n+1}(A)$ sont donn\'ees 
pour $a_0, a_1, \ldots, a_n \in A$ et $i=0,\ldots, n$ par
$$d_i(a_0\ot a_1\ot \cdots \ot a_n) = 
\left\{
\matrix{
a_0\ot \cdots \ot a_ia_{i+1}\ot \cdots \ot a_n & 
\quad \hbox{si}\; 0\leq i \leq n-1, \hfill\cr
\noalign{\medskip}
a_n a_0\ot a_1\ot \cdots \ot a_{n-1} &
\quad \hbox{si}\; i = n, \hfill\cr
}
\right.
$$
et
$$s_i(a_0\ot a_1\ot \cdots \ot a_n) = 
a_0\ot \cdots \ot a_i\ot 1\ot a_{i+1}\ot \cdots \ot a_n .
$$

\medskip
\noindent
{\sc 5.9.\ Proposition.} 
{\it Soit $B = A\ot A^{\op}$ l'alg\`ebre enveloppante de l'alg\`ebre
associative uni\-f\`ere~$A$. 
Posons $P_n = A^{\ot (n+2)}$
si $n\geq 0$ et $P_n = 0$ si $n<0$.
On d\'efinit $\eps : P_0 \to A$ et $d':P_n\to P_{n-1}$ par
$$\eps(a_0\ot a_1) = a_0 a_1 \et
d'(a_0\ot a_1\ot \cdots \ot a_{n+1}) = 
\sum_{i=0}^{n}\, q^i\, a_0\ot \cdots \ot a_ia_{i+1}\ot \cdots \ot a_{n+1}$$
o\`u $a_0, a_1, \ldots, a_n \in A$.
Alors $(P,d', \eps : P_0\to A)$ est une $N$-r\'esolution
de~$A$ pour la classe projective $k$-scind\'ee~$(\cP_k,\cS_k)$
de~${}_B\mod$.
}

\medskip
\Dem
Les $B$-modules $P_n$ sont des modules projectifs relativement aux
morphismes $k$-scind\'es  puisqu'ils sont de la forme 
$A\ot A^{\ot n}\ot A \cong B\ot A^{\ot n}$. 
L'acyclicit\'e du $N$-complexe $(P,d',\eps:P_0\to A)$ est
une cons\'equence de la proposition~5.2.
\hfill\cqfd

\medskip 
\noindent
{\sc 5.10.\ D\'emonstration du Th\'eor\`eme~1.}---
Consid\'erons la $N$-r\'esolution $\cP_k$-projective $(P,d', \eps : P_0\to A)$
de $B$-modules fournie par la proposition~5.9.
Un calcul facile montre que le $N$-complexe $(P\ot _{B} A, d'\ot \id)$
est isomorphe au $N$-complexe~(0.1) $(C_*(A),b)$ de l'introduction.
Par cons\'equent, son homologie co\"\i ncide avec les avatars de l'homologie
de Hochschild d\'efinis par~(0.2)~:
$${}_pH_*(P\ot _{B} A, d'\ot \id) \cong {}_pHH_*(A).$$
Or, d'apr\`es la d\'efinition~3.2 appliqu\'ee \`a la m\^eme
$N$-r\'esolution $\cP_k$-projective, on a
$${}_pH_*(P\ot _{B} A, d'\ot \id) \cong {}_p\Tor_*^B(A,A),$$
d'o\`u nous tirons l'isomorphisme~(0.3) de l'introduction.
On termine la d\'emonstration du th\'eor\`eme~1 en utilisant le
corollaire~3.3 et en rappelant que l'homologie de Hochschild
usuelle (c'est-\`a-dire celle du complexe de Hochschild standard)
est isomorphe \`a~$\Tor_*^B(A,A)$.
\hfill\cqfd

\medskip 
\noindent
{\sc 5.11.\ Remarque.}---
Les m\^emes m\'ethodes permettent de donner une version cohomo-logique
du th\'eor\`eme~1. Nous laissons au lecteur le soin de l'\'enoncer.

\goodbreak

\vskip 25pt
\vfill\eject


\centerline{\bf R\'ef\'erences}

\vskip 20pt
\vskip 10pt

\noindent
[Dub] {M.~Dubois-Violette}, 
{\it Generalized differential spaces with $d^N=0$ 
and the $q$-differential calculus}, 
preprint q-alg/9609012.
\vskip 3pt 

\noindent
[DK1] {M.~Dubois-Violette, R.~Kerner}, 
{\it Universal $q$-differential calculus 
and $q$-analog of homological algebra}, 
preprint q-alg/9608026.
\vskip 3pt 

\noindent
[DK2] {M.~Dubois-Violette, R.~Kerner}, 
{\it Universal $Z_N$-graded differential calculus}, 
preprint.
\vskip 3pt

\noindent
[EM] {S.~Eilenberg, J. C. Moore}, 
{\it Foundations of relative homological algebra}, 
Memoirs of the Amer.\ Math.\ Soc.\ 55, Providence, Rhode Island, 1965.
\vskip 3pt

\noindent
[Kap] {M.~M.\ Kapranov}, 
{\it On the $q$-analog of homological algebra}, 
preprint q-alg/9611005.
\vskip 3pt

\noindent
[May] {W.~Mayer}, 
{\it A new homology theory, I, II}, 
Annals of Math.\ 43 (1942), 370--380 et 594--605.
\vskip 3pt

\noindent
[Sar1] {K. S. Sarkaria}, 
{\it Combinatorial methods in topology}, 
Notes of Chandigarh Topology Seminar (1994--95).
\vskip 3pt

\noindent
[Sar2] {K. S. Sarkaria}, 
{\it Some simplicial (co)homologies}, 
preprint IHES/M/95/83.
\vskip 3pt

\noindent
[Spa] {E. H. Spanier}, 
{\it The Mayer homology theory}, 
Bull.\ Amer.\ Math.\ Soc.\ 55 (1949), 102--112.
\vskip 3pt

\vskip 10pt
\vskip 20pt

\line{Institut de Recherche Math\'ematique Avanc\'ee \hfill}
\line{Universit\'e  Louis Pasteur - C.N.R.S.\hfill}
\line{7 rue Ren\'e Descartes \hfill}
\line{67084  Strasbourg Cedex, France \hfill}
\line{E-mail~:  {\tt kassel@math.u-strasbg.fr, wambst@math.u-strasbg.fr}\hfill}
\line{Fax~: +33 (0)3 88 61 90 69 \hfill}

\vskip 30pt
\noindent
(29 avril 1997)
\bye